\documentclass[11pt]{amsart}

\usepackage[utf8x]{inputenc} 
\usepackage[a4paper]{geometry} 
\usepackage{enumerate}  
\usepackage{fancyvrb}    
\usepackage{color,colortbl}
\usepackage[dvipsnames]{xcolor}
\usepackage{hyperref}  
\usepackage{verbatim}
                 \hypersetup{ pdfborder={0 0 0}, 
                              colorlinks=true, 
                              linkcolor=red, 
                              urlcolor=red,
                              citecolor=red,
                              linktoc=page,
                              pdfauthor={Giovanni Staglian{\`o}}, 
                              pdftitle={A \emph{Macaulay2} package for computations with rational maps}
                            }                         
\usepackage{mathptmx}

\theoremstyle{plain} 
\newtheorem{proposition}{Proposition}[section] 
\theoremstyle{definition} 
\newtheorem{definition}[proposition]{Definition} 
\newtheorem{example}[proposition]{Example} 
\newtheorem{algorithm}[proposition]{Algorithm}
\theoremstyle{remark} 
\newtheorem{remark}[proposition]{Remark}

\newcommand{\B}{{\mathfrak{B}}}

\newcommand{\PP}{{\mathbb{P}}}

\numberwithin{equation}{section}

\title[A Macaulay2 package for computations with rational maps]{A Macaulay2 package for computations with rational maps} 
\author[G. Staglian\`o]{Giovanni Staglian\`o} 
\date{\today} 
\address{Dipartimento di Ingegneria Industriale e Scienze Matematiche, Universit\`a Politecnica delle Marche, via Brecce Bianche, I-60131 Ancona, Italy}
\email{\href{mailto:giovannistagliano@gmail.com}{giovannistagliano@gmail.com}} 

\begin{document}

\begin{abstract} 
The \emph{Macaulay2} package \emph{Cremona} performs some 
computations on rational and birational maps between 
irreducible projective varieties. 
For instance, it provides methods to compute degrees 
and projective degrees of rational maps without any theoretical limitation,
from which is derived a general method 
to compute the push-forward to projective space of Segre classes.
Moreover, the computations can be done both deterministically and probabilistically. 
We give here a brief description of the 
methods and algorithms implemented. 
\end{abstract}

\maketitle

\section*{Introduction} 
In this note we describe the computational package  
\href{https://github.com/Macaulay2/M2/blob/master/M2/Macaulay2/packages/Cremona.m2}{\emph{Cremona}},  
included with \emph{Macaulay2} since version 1.9 \cite{macaulay2}.
A first rudimentary version of this package 
has been already used 
  in an essential way in \cite{note3}  
(it was originally named \href{http://goo.gl/eT4rCR}{bir.m2}),
and recent applications can be found in \cite{note4Cubics,russo-stagliano}.
Here we describe version 4.2.1 of the package,
 included with \emph{Macaulay2} version 1.11.
  
\emph{Cremona}
 performs 
 computations on rational and 
birational maps between absolutely 
irreducible projective 
varieties  over a field $\mathbb{K}$. Among other things, it provides
general methods
to compute 
projective degrees 
of rational maps, 
from which, as is well-known (see Proposition~\ref{segreprop}), 
one can interpret them as methods to compute 
  the push-forward to projective space of Segre classes.
  The algorithms  are 
  naively derived 
  from the mathematical definitions, with the advantages of being
  obvious, quite general and easily implemented. 
  Moreover, all the methods (where this may make sense) 
are available both in a probabilistic version  
and in a deterministic version, 
and one can switch from one to the other with a boolean option named \texttt{MathMode}.

In Section \ref{sec1}, we will describe the main methods provided by the package
and the algorithms implemented. 
Most of these have already been described in \cite[Section~2]{note3}, 
but here we will consider a more general setting. For instance, 
Algorithm \ref{compKer} for computing homogeneous components of kernels of homogeneous ring maps  
was presented in \cite[Algorithm~2.5]{note3} requiring that the map was between polynomial rings. 
In section \ref{sec2}, we will show how these methods work in  some particular examples,   
 concluding with an experimental comparison of 
  the running times  of one of these  methods 
  with the corresponding ones proposed 
  in \cite{Helmer2016120} and \cite{caroyharris} (see also \cite{Jost2015}).
    For further technical details we refer to
   the documentation of the package,  
which can be shown using the command \texttt{viewHelp Cremona}.

  We mention that the package 
\href{https://github.com/Macaulay2/M2/blob/master/M2/Macaulay2/packages/RationalMaps.m2}{\emph{RationalMaps}},
by K. Schwede,
D. Smolkin, 
S. H. Hassanzadeh, and
C. J. Bott, 
is another package included with \emph{Macaulay2} for working with rational maps. 
It mainly focuses on providing a general method for inverting birational maps, which in some cases
turns out to be competitive with the corresponding method of \emph{Cremona}.

\section{Description of the main methods}\label{sec1}

Throughout, we shall use the following notation. 
Let $\mathbb{K}$ denote a field; in practical, it can be for instance 
$\mathbb{Q}$, a finite field, or a fraction field of a polynomial ring over these.
Let $\phi:X\dashrightarrow Y$ be a rational map 
from a subvariety $X=V(I)\subseteq\mathbb{P}^n=\mathrm{Proj}(\mathbb{K}[x_0,\ldots,x_n])$ 
to a subvariety $Y=V(J)\subseteq\mathbb{P}^m=\mathrm{Proj}(\mathbb{K}[y_0,\ldots,y_m])$, 
which can be represented, although not uniquely, 
by a homogeneous ring map 
$\varphi:\mathbb{K}[y_0,\ldots,y_m]/J\rightarrow\mathbb{K}[x_0,\ldots,x_n]/I$ 
of quotients of polynomial rings by homogeneous ideals. 
Sometimes we will denote by $F_0,\ldots,F_m\in \mathbb{K}[x_0,\ldots,x_n]$ 
homogeneous forms of the same degree
such that $\bar{F_i} := F_i + I = \varphi(y_i)$, for $i=0,\ldots,m$.
The common degree of these elements will be denoted by $\delta$.

\subsection{From algebraic geometry to computational algebra}
For each homogeneous ideal $\mathfrak{a}\subseteq \mathbb{K}[x_0,\ldots,x_n]/I$ (resp. 
 $\mathfrak{b}\subseteq \mathbb{K}[y_0,\ldots,y_m]/J$), we have 
a closed subscheme $V(\mathfrak{a})\subseteq X$ (resp.
$V(\mathfrak{b})\subseteq Y$), and
the following 
basic
formulae hold:\footnote{By a little abuse of notation, we here are considering $\phi$ as a morphism defined  on the open set $X\setminus V(\bar{F_0},\ldots,\bar{F_m})$.}
\begin{equation}\label{formule}
    \overline{\phi\left( V(\mathfrak{a})\right)} = V(\varphi^{-1}(\mathfrak{a}))\quad \mbox{ and }\quad 
    \overline{\phi^{-1}\left( V(\mathfrak{b})\right)} = V\left(\left(\varphi(\mathfrak{b})\right) : {\left(\varphi(y_0),\ldots,\varphi(y_m)\right)}^{\infty}\right) .
\end{equation}
In particular, the (closure 
of the) image of $\phi$ is defined by the kernel of $\varphi$.
Several issues about rational maps lead naturally to
an examination of the left-hand sides of \eqref{formule}, and the right-hand sides of \eqref{formule}
can be determined  using Gr\"{o}bner basis techniques, 
whenever $\mathfrak{a}$ and $\mathfrak{b}$ are explicitly given. 
Furthermore, \emph{Macaulay2} provides 
useful commands as \texttt{preimage}, \texttt{kernel} and \texttt{saturate}, 
so that the required users' skills are quite low. 
The aim of the package \emph{Cremona} is to provide further tools. 

\subsection{Computing projective degrees} 
The projective degrees are the most basic invariants of a rational map. 
Many others can be derived from them, as for instance the dimension and the degree of the base locus. For more details on the subject, see \cite[Example~19.4, p.~240]{harris-firstcourse}.
\begin{definition}[\emph{Projective degrees}, \cite{harris-firstcourse}] \label{defproj} \hspace{1pt}
 \begin{enumerate}
  \item\label{multid} The projective degrees $d_0(\phi),d_1(\phi),\ldots,d_{\dim X}(\phi)$ of the map $\phi$ are 
  defined as the components of 
  the multidegree of the closure of the graph $\Gamma_{\phi}\subset \mathbb{P}^n\times\mathbb{P}^m$; 
  equivalently,
  \item\label{harrisdeg} the $i$-th projective degree $d_i(\phi)$ can be defined in terms of dimension 
and degree of the closure of $\phi^{-1}(L)$, 
 where $L$ is a general  $(m-\dim X + i)$-dimensional linear 
 subspace of $\mathbb{P}^m$; precisely, $d_i(\phi)=\deg \overline{\phi^{-1}(L)}$ 
 if $\dim \overline{\phi^{-1}(L)} = i$, and $d_i(\phi)=0$ otherwise.
 \end{enumerate}
\end{definition}
In common computer algebra systems such as \emph{Macaulay2}, 
it is easy to translate Definition~\ref{defproj} into code.
 We now describe more in details
how this can be done. 
All of this is implemented in the method \texttt{projectiveDegrees};
see Example \ref{cubo} for an example using it.
\subsubsection{Deterministic approach} 
Taking into account Definition \ref{defproj}(\ref{multid}),
a bihomogeneous ideal for $\Gamma_{\phi}$ in
$\mathbb{K}[x_0,\ldots,x_n,y_0,\ldots,y_m]$
can be, for instance, obtained as
\begin{equation}
\left(I + \left(\left\{y_i\,F_j-y_j\,F_i,\  0\leq i,j\leq m\right\}\right)\right) : {(F_0,\ldots,F_m)}^{\infty}.
\end{equation}
Therefore its  multidegree   
can be computed in \emph{Macaulay2} with the command \texttt{multidegree},
which implements an algorithm according to \cite[p.~165]{millersturmfels}.

\subsubsection{Probabilistic approach} (See also \cite[Remark~2.4]{note3}.) 
Taking into account Definition \ref{defproj}(\ref{harrisdeg}),
if 
$L$ is defined by an ideal $I_L$, the second formula of \eqref{formule} says us 
 that $\overline{\phi^{-1}(L)}$ 
is defined by the saturation
of the ideal $\left(\varphi(I_L)\right)$ by $(\bar{F_0},\ldots,\bar{F_m})$ 
in the ring $\mathbb{K}[x_0,\ldots,x_n]/I$. So replacing the word 
\emph{general} with \emph{random} in the definition,  we get a probabilistic algorithm that computes
all the projective degrees. 
Moreover, we can considerably speed up this algorithm by taking into account 
two 
remarks: 
firstly, the saturation $\varphi(I_L):{(\bar{F_0},\ldots,\bar{F_m})}^{\infty}$ 
is the same as $\varphi(I_L):{(\lambda_0\,\bar{F_0}+\cdots+\lambda_m\,\bar{F_m})}^{\infty}$, 
where $\lambda_0,\ldots,\lambda_m\in\mathbb{K}$ are general  
scalars; 
secondly,
the $i$-th projective degree of $\phi$ coincides with the $(i-1)$-th projective degree 
of the restriction of $\phi$ to a general 
hyperplane section of $X$.

\subsubsection{An alternative deterministic approach} 
Replacing 
 the word 
\emph{general} with \emph{symbolic} 
in  Definition \ref{defproj}(\ref{harrisdeg})
gives us a deterministic algorithm for computing 
 projective degrees. 
For instance, 
in the case in which $\phi:\PP^n\dashrightarrow \PP^n$ is a dominant rational map, 
extending $\mathbb{K}$ to the fractional field of a polynomial ring
$\mathbb{K}[a_0,\ldots,a_n]$, we have that 
$d_0(\phi)$ is the degree of the fibre of $\phi$  
 at the \emph{symbolic point} $[a_0,\ldots,a_n]$.

\subsection{Some applications using projective degrees}
\subsubsection{The degree of a rational map} 
The degree of the  map $\phi:X\dashrightarrow Y$ 
is the number of isolated points  in the inverse image of a general point of $\overline{\phi(X)}$ 
over the algebraic closure of $\mathbb{K}$. This is the same as the 
ratio of $d_0(\phi)$ and  $\deg \overline{\phi(X)}$, and  
thus it can be explicitly computed.
Let us note, however, that in several cases
we do not need to compute the kernel of $\varphi$. 
For instance, if $X$ is a projective space,
we are able to pick up 
an abundance of rational points
 of $\overline{\varphi(X)}$ and 
then we apply the second formula of \eqref{formule}.
Another special case is when $d_0(\phi)$ is a prime number:
here we have only to establish if the image of $\phi$ is a linear subspace 
(e.g. applying Algorithm \ref{compKer} with $d=1$).
The method provided by \emph{Cremona} 
for this computation is named   \texttt{degreeOfRationalMap}.\footnote{Notice that,
in general, 
if the result of the probabilistic algorithm for \texttt{degreeOfRationalMap} 
is wrong, it can be either to small or to large. 
However, as a consequence of \cite[Chapter~III, Exercise~10.9]{hartshorne-ag}, it should always provide a lower bound when the map is dominant between smooth varieties.} 

Methods related to this are  \texttt{isBirational} and \texttt{isDominant}  
with obvious meaning. The latter does not compute the kernel of $\varphi$, but 
it uses 
an algorithm that looks for  $d_r(\phi)$, where $r =\dim X - \dim Y$. 
 More precisely, the algorithm is based on the following fact: 
    let $Z\subset Y$ be a 
          random zero-dimensional linear section of $Y$; if
          $\dim \overline{\phi^{-1}(Z)} = \dim X - \dim Y \geq 0$, 
          then $\phi$ is certainly dominant, otherwise it is probably not dominant 
          (see \cite[Chapter~I, \S~8]{mumford} or \cite[Chapter~II, Exercise~3.22]{hartshorne-ag}). 
          When this last case occurs, 
     it is generally easy 
       to find a nonzero element in the kernel of $\varphi$,
       and so this method turns out to be very 
       effective even
       in its deterministic version (see Example \ref{cre20}).
  
\subsubsection{The Segre class} 
It is well-known that one can deduce an algorithm computing 
the push-forward to projective space of Segre classes 
from an algorithm computing projective degrees 
of rational maps between projective varieties and vice versa. 
Indeed, with our notation, we have the following:
\begin{proposition}[Proposition~4.4 in \cite{fulton-intersection}; 
             see also Subsection~2.3 in \cite{dolgachev-cremonas} and Section~3 in \cite{aluffi}] 
  \label{segreprop}
  Let $\B\subset X$ be the subscheme defined by $\bar{F_0},\ldots,\bar{F_m}$ 
  and 
  let 
  $\nu:X\hookrightarrow \PP^n$ be the inclusion. If $H$ denotes the hyperplane class of $\PP^n$ and
 $r = \dim X$,
 then 
 the push-forward $\nu_{\ast}(s(\B,X))$ of the Segre class of $\B$ in $X$ is 
 \begin{equation}\label{segre}
\nu_{\ast}(s(\B,X))  =   \sum_{k=0}^{\dim \B}\left( (-1)^{r -k-1}\,\sum_{i=0}^{r -k} (-1)^i \binom{r -k}{i}\, \delta^{r -k-i}\, d_{r-i}(\varphi)\right)\, H^{n-k}  .  
 \end{equation} 
\end{proposition}
The 
general method \texttt{SegreClass}, provided by \emph{Cremona}
for computing the push-forward to projective space of Segre classes,
does basically nothing more than apply \eqref{segre};
 see Example \ref{segreDet} for an example using this method.
 Furthermore, 
applying one of the main results in \cite{aluffi}, 
a method is derived to compute
the push-forward to projective space of the
  Chern-Schwartz-MacPherson class $c_{SM}(W)$ of the support of a projective scheme $W$;
  recall that the component of dimension $0$ of $c_{SM}(W)$ is the topological Euler 
  characteristic of the support of $W$.
    

\subsection{Computing homogeneous components of kernels}\label{compkernels}
To compute, using \emph{Macaulay2}, the homogeneous component of degree $d$ of the kernel 
of 
\texttt{phi=}$\varphi$, one can perform the command 
\texttt{ideal image basis(d,kernel phi)}.
This is equivalent to the command
\texttt{kernel(phi,d)} provided by \texttt{Cremona},
but the latter uses the following obvious algorithm.
\begin{algorithm}\label{compKer}    \hspace{1pt}
 \begin{description} 
  \item[\texttt{Input}] the ring map $\varphi$ and an integer $d$.
 \item[\texttt{Output}] homogeneous component of degree $d$ of the kernel of $\varphi$.
 \end{description}
\begin{itemize} 
  \item Find vector space bases $G_0,\ldots,G_r$
  of $\left(\mathbb{K}[y_0,\ldots,y_m]/J\right)_d$ and $H_0,\ldots,H_s$ of $I_{d\,\delta}$,
  where subscripts stand for homogeneous components; 
  \item take generic linear combinations 
$\mathbf{G}=\sum_{i=0}^{r} a_i\,G_i$ and 
 $\mathbf{H}=\sum_{j=0}^{s}b_j\,H_j$, and      
  find a basis of solutions for the
  homogeneous linear system obtained by imposing that
  the polynomial 
  $\mathbf{G}(F_0,\ldots,F_m) - \mathbf{H}\in\mathbb{K}[a_0,\ldots,a_r,b_0,\ldots,b_s][x_0,\ldots,x_n]$  
  vanishes identically;
 \item   
 for each vector $(\hat{a_0},\ldots,\hat{a_r},\hat{b_0},\ldots,\hat{b_s})\in\mathbb{K}^{r+s+2}$ 
 obtained in the previous step, replace in
        $\mathbf{G}$ the coefficients $a_0,\ldots,a_r$ with $\hat{a_0},\ldots,\hat{a_r}$;
        return 
        all these elements.
 \end{itemize}
\end{algorithm}
For small values of $d$, applying Algorithm \ref{compKer}  may turn out to be much faster than
computing a list of generators of the kernel of the map;
see for instance Example \ref{cre20} below.

\subsection{Inverting birational maps}
General algorithms for inverting birational maps are known. 
One of them is implemented in the package 
 \href{https://github.com/Macaulay2/M2/blob/master/M2/Macaulay2/packages/Parametrization.m2}{\emph{Parametrization}} by J. Boehm, and the method \texttt{inverseMap} of \emph{Cremona} 
 uses the same one for the general case as well. However,
when the source $X$ of the rational map $\phi$ is a projective space 
and a further technical condition is satisfied, 
then it uses the following powerful algorithm. 
\begin{algorithm}[\cite{russo-simis}; see also \cite{Simis2004162}]\label{inverseRS} 
 \hspace{1pt}
 \begin{description} 
  \item[\texttt{Input}] the ring map $\varphi$ 
  (assuming that $\phi$ is birational and further conditions are satisfied).
 \item[\texttt{Output}] a ring map representing the inverse map of $\phi$. 
 \end{description}
\begin{itemize} 
 \item   Find
 generators $\{(L_{0,j},\ldots,L_{m,j})\}_{j=1,\ldots,q}$ for  
  the module of linear syzygies of 
  $(F_0,\ldots,F_m)$;
  \item compute the Jacobian matrix $\Theta$ 
of the bihomogeneus forms $\{\sum_{i=0}^m y_i\,L_{i,j} \}_{j=1,\ldots,q}$ 
with  respect to the variables $x_0,\ldots,x_n$
and consider 
the map of graded free modules 
$(\Theta\ \mathrm{mod}\ J):
\left(k[y_0,\ldots,y_m]/J\right)^{n+1}\rightarrow 
\left(k[y_0,\ldots,y_m]/J\right)^{q}$;
\item return the map defined by 
a generator $G=(G_0,\ldots,G_n)$ for the kernel of $(\Theta\ \mathrm{mod}\ J)$.
\end{itemize}
\end{algorithm}
\begin{remark}
One of the main feature of
the package 
\href{https://github.com/Macaulay2/M2/blob/master/M2/Macaulay2/packages/RationalMaps.m2}{\emph{RationalMaps}},
by K. Schwede,
D. Smolkin, 
S. H. Hassanzadeh, and
C. J. Bott, 
is  a method for inverting birational maps, which, 
 in the case when Algorithm~\ref{inverseRS}  
 does not apply, 
  appears to be quite competitive
 with
 the method \texttt{inverseMap} of \emph{Cremona}.
\end{remark}

\subsubsection{Heuristic approach} 
The method \texttt{approximateInverseMap} 
provides a heuristic approach to compute 
the inverse of a birational map modulo a change of coordinate.
The idea of the algorithm is to try to construct 
the base locus of the inverse by looking for the images 
of general linear sections.
Consider, for simplicity, the case in which $\phi:\PP^n\dashrightarrow{\PP^n}'$
is a Cremona transformation.
Then, by taking the images of $n+1$ 
general hyperplanes in $\PP^n$, 
we form a linear system of hypersurfaces in ${\PP^n}'$ of degree $d_1(\phi)$ 
which defines a rational map $\psi:{\PP^n}'\dashrightarrow\PP^n$ such that 
 $\psi\circ \phi$ 
 is a (linear) isomorphism; i.e. we find an \emph{approximation} of $\phi^{-1}$.
Next, 
we can fix the \emph{error of the approximation} 
by observing that we have
$\phi^{-1} = \left(\psi\circ \phi\right)^{-1}\circ \psi$. 
It is surprising that this method turns to be effective in 
 examples where other deterministic algorithms seem to run endlessly; see for instance 
 Example \ref{cre20} below.
\section{Examples}\label{sec2}
In this section, we show how the methods described in Section \ref{sec1} 
can be applied in some particular examples. 
We point out that the package \emph{Cremona} provides the data type \texttt{RationalMap}, but here 
we will use the more familiar 
     type \texttt{RingMap}. For brevity, we will omit irrelevant output lines.
We start with an example reviewing 
the construction given in \cite{note3} 
 of a quadro-quadric Cremona transformation of $\PP^{20}$.
 
 \begin{example}\label{cre20}
 The code below constructs a ring map  \texttt{psi} 
 representing 
 a rational map $\psi:\PP^{16}\dashrightarrow\PP^{20}$.%
\footnote{Precisely, the algorithm for constructing $\psi$ is as follows:
 take $E\subset\PP^7$ to be a 
$3$-dimensional Edge variety of degree $7$, namely, 
the residual intersection of $\PP^1\times\PP^3\subset\PP^7$ with
a general quadric in $\PP^7$ containing one of the $\PP^3$'s of the 
rulings of $\PP^1\times\PP^3\subset\PP^7$;
next, see $E\subset\PP^7$ embedded in a hyperplane of $\mathbb{P}^8$ 
and take the birational map 
$\phi:\PP^8\dashrightarrow \PP^{16}$
defined by the quadrics of $\PP^8$ containing $E$;
 take $\psi:\PP^{16}\dashrightarrow\PP^{20}$ 
to be the map defined by the quadrics of $\PP^{16}$ containing 
the image of $\phi$.} 
For this first part, 
we use the package  \emph{Cremona}
 only to shorten the code.
{\footnotesize
\begin{Verbatim}[commandchars=&\{\}]
Macaulay2, version 1.11
with packages: ConwayPolynomials, Elimination, IntegralClosure, InverseSystems, 
               LLLBases, PrimaryDecomposition, ReesAlgebra, TangentCone
i1 : loadPackage "&colore{blue}{Cremona}";
i2 : K = ZZ/70001;
i3 : P8 = K[t_0..t_8]; E = saturate(minors(2,genericMatrix(P8,4,2))+sum(
     (ideal(t_0..t_7)*ideal(t_0..t_3))_*,u->random(K)*u),ideal(t_0..t_3)) + t_8;
i5 : psi = &colore{blue}{toMap kernel(toMap}(E,2),2);
\end{Verbatim}
} \noindent 
Up to this point, the computation was standard. But now
we want to determine the homogeneous ideal of 
$Z:=\overline{\psi(\PP^{16})}\subset\PP^{20}$, 
which turns out to be generated by quadrics.
Computing this using \texttt{kernel psi}
seems an impossible task, but it is elementary using
  \texttt{kernel(psi,2)}.
So we can consider $\psi$ as a dominant rational map 
$\psi:\PP^{16}\dashrightarrow Z\subset \PP^{20}$.
{\footnotesize
\begin{Verbatim}[commandchars=&\{\}]
i6 : time Z = &colore{blue}{kernel}(psi,2);
     &colore{Sepia}{-- used 2.84998 seconds}
i7 : psi = &colore{blue}{toMap}(psi,&colore{blue}{Dominant}=>Z);
\end{Verbatim}
} \noindent 
The map $\psi$ turns out to be not only dominant but birational. 
{\footnotesize
\begin{Verbatim}[commandchars=&\{\}]
i8 : time &colore{blue}{degreeOfRationalMap} psi
     &colore{Sepia}{-- used 2.11216 seconds}
o8 = 1
\end{Verbatim}
} \noindent 
We now want to compute the inverse of $\psi$.
This is a case where \texttt{inverseMap} can apply
  Algorithm \ref{inverseRS},
but the running time is several hours.
We can perform this computation in seconds by using \texttt{approximateInverseMap}.
{\footnotesize 
\begin{Verbatim}[commandchars=&\{\}] 
i9 : time psi' = &colore{blue}{approximateInverseMap}(psi,&colore{blue}{CodimBsInv}=>10,&colore{blue}{MathMode}=>true);
     &colore{Sepia}{-- used 15.9724 seconds}
\end{Verbatim}
} \noindent 
A Cremona transformation $\omega$ of $\PP^{20}$ is then obtained combining 
$\psi^{-1}$ and $Z$ as follows.
{\footnotesize
\begin{Verbatim}[commandchars=&\{\}] 
i10 : omega = &colore{blue}{toMap}(lift(matrix psi',ring Z)|gens Z);
\end{Verbatim}
} \noindent 
Even checking just the dominance of $\omega$, by computing \texttt{kernel omega}, seems an
  impossible task, but it can be done quickly with \texttt{isDominant}.
{\footnotesize
\begin{Verbatim}[commandchars=&\{\}] 
i11 : time &colore{blue}{isDominant}(omega,&colore{blue}{MathMode}=>true)
     &colore{Sepia}{-- used 0.100369 seconds}
o11 = true
\end{Verbatim}
} \noindent 
We now check that our map is birational and compute its inverse
using Algorithm \ref{inverseRS}.
{\footnotesize
\begin{Verbatim}[commandchars=&\{\}] 
i12 : time &colore{blue}{isBirational} omega
     &colore{Sepia}{-- used 0.0366468 seconds}
o12 = true
i13 : time &colore{blue}{inverseMap} omega;
     &colore{Sepia}{-- used 0.0717518 seconds}
\end{Verbatim}
} \noindent 
\end{example}
\begin{example}\label{cubo} 
In this example, we use the probabilistic versions of some methods.
Take $M$ to be a generic $3 \times 5$ matrix of linear forms on $\PP^6$, 
and let $\phi:\PP^6\dashrightarrow\mathbb{G}(2,4)\subset\mathbb{P}^9$ be 
the rational map defined by the $3\times 3$ minors of $M$ (its base locus 
is a smooth threefold scroll 
over a plane).
{\footnotesize
\begin{Verbatim}[commandchars=&\{\}] 
i14 : P6 = K[x_0..x_6]; M = matrix pack(5,for i from 1 to 15 list random(1,P6));
i16 : phi = &colore{blue}{toMap}(minors(3,M),&colore{blue}{Dominant}=>2);
\end{Verbatim}
} \noindent 
We check that the map is birational and compute its inverse.
{\footnotesize
\begin{Verbatim}[commandchars=&\{\}]  
i17 : time &colore{blue}{isBirational} phi
     &colore{Sepia}{-- used 0.217607 seconds}
o17 = true
i18 : time psi = &colore{blue}{inverseMap} phi;
    &colore{Sepia}{ -- used 1.39511 seconds}
\end{Verbatim}
} \noindent 
Now we compute the multidegrees of $\phi$ and $\phi^{-1}$. 
{\footnotesize
\begin{Verbatim}[commandchars=&\[\]]
i19 : time (&colore[blue][projectiveDegrees] phi, &colore[blue][projectiveDegrees] psi)
     &colore[Sepia][-- used 1.37582 seconds]
o19 = ({1, 3, 9, 17, 21, 15, 5}, {5, 15, 21, 17, 9, 3, 1})
\end{Verbatim}
} \noindent 
We also compute the push-forward to $\PP^6$ (resp. $\PP^9$)
of the Segre class of the base locus of $\phi$ (resp. $\phi^{-1}$) in
$\PP^6$ (resp. in $\mathbb{G}(2,4)$). As usual, $H$ denotes the hyperplane class.
{\footnotesize
\begin{Verbatim}[commandchars=&\{\}] 
i20 : time (&colore{blue}{SegreClass} phi, &colore{blue}{SegreClass} psi)
     &colore{Sepia}{-- used 1.43359 seconds}
             6       5      4      3      9       8       7      6      5
o20 = (- 680H  + 228H  - 60H  + 10H , 728H  - 588H  + 276H  - 98H  + 24H )
\end{Verbatim}
} \noindent 
\end{example} 
\begin{example}\label{segreDet} 
 In this example, we use the deterministic version of the method \texttt{SegreClass}.
 We take $Y\subset\PP^{11}$ to be the dual quartic hypersurface of 
$\PP^1\times Q^4\subset{\PP^{11}}^{\ast}$, 
where $Q^4\subset\PP^5$ is a smooth quadric hypersurface,
and take $X\subset Y$ to be the singular locus of $Y$. 
We then compute the push-forward to the Chow ring of $\PP^{11}$
of the Segre class both of $X$ in $Y$ and of $X$ in $\PP^{11}$ working over 
the  Galois field $\mathrm{GF}(331^2)$.
{\footnotesize
\begin{Verbatim}[commandchars=&\{\}]  
i21 : P11 = GF(331^2)[x_0..x_11];
i22 : Y = ideal sum(first entries gens minors(2,genericMatrix(P11,6,2)),t->t^2);
i23 : X = sub(ideal jacobian Y,P11/Y);
i24 : time &colore{blue}{SegreClass}(X, &colore{blue}{MathMode}=>true)  &colore{Sepia}{-- push-forward of s(X,Y)}
     &colore{Sepia}{-- used 0.789986 seconds}
             11          10         9        8        7       6       5      4      3
o24 = 507384H   - 137052H   + 35532H  - 9018H  + 2340H  - 658H  + 204H  - 64H  + 16H
i25 : time &colore{blue}{SegreClass}(lift(X,P11), &colore{blue}{MathMode}=>true)  &colore{Sepia}{-- push-forward of s(X,P^11)}
    &colore{Sepia}{ -- used 0.846234 seconds}
             11          10         9        8        7       6       5      4     3
o25 = 313568H   - 101712H   + 30636H  - 8866H  + 2532H  - 720H  + 198H  - 48H  + 8H
\end{Verbatim}
} \noindent 
\end{example}
\begin{example} 
Here
we experimentally measure the probability of obtaining an incorrect answer
using the probabilistic version of the method \texttt{projectiveDegrees} 
with a simple example of a birational map 
$\phi:\mathbb{G}(1,3)\dashrightarrow\PP^4$ defined over $\mathbb{K}$.
We define a procedure which computes this probability as a function of the field $\mathbb{K}$. 
In Table~\ref{table: probability}, we report the results obtained by running the procedure with various fields.

{\footnotesize
\begin{Verbatim}[commandchars=&?$]  
i26 : p = (K) -> (
         x := local x; R := K[x_0..x_4];
         phi := &colore?blue$?inverseMap toMap$(minors(2,matrix{{x_0..x_3},{x_1..x_4}}),&colore?blue$?Dominant$=>2);
         m := &colore?blue$?projectiveDegrees$(phi,&colore?blue$?MathMode$=>true);
         0.1 * # select(1000,i -> &colore?blue$?projectiveDegrees$ phi != m));
\end{Verbatim}     
} \noindent 

 \begin{table}[htbp]
\centering 
\renewcommand{\arraystretch}{1.1}
\begin{tabular}{|l|rrrrr|}
 \hline 
  Field:  & $\mathbb{Q}$ & $\mathbb{Z}/70001$ & $\mathrm{GF}(3^8)$ & $\mathbb{Z}/101$ & $\mathbb{Z}/31$ \\
 Probability:  &  0.0\% &  0.0\% &  0.2\% &  7.4\% &  25.3\%\\ 
\hline 
\end{tabular}
 \caption{Incorrect outputs of a probabilistic method.}
\label{table: probability} 
\end{table}

\end{example}

\begin{example}\label{aluffihelmer}
In this last example, 
 we deal  with an experimental comparison 
  of the method \texttt{SegreClass} of \emph{Cremona}
  and the corresponding ones of other \emph{Macaulay2} packages.
Precisely, 
we want to compare the method 
\texttt{SegreClass} against 
the corresponding methods 
 of the two following packages:
 \href{https://github.com/Macaulay2/M2/blob/master/M2/Macaulay2/packages/CharacteristicClasses.m2}{\emph{CharacteristicClasses}} version 2.0, 
 by M. Helmer and C. Jost (see \cite{Helmer2016120,Jost2015}), which provides a probabilistic method; and
 \href{https://github.com/coreysharris/FMPIntersectionTheory-M2}{\emph{FMPIntersectionTheory}} version 0.1,
 by C. Harris (see \cite{caroyharris}), which provides a deterministic 
 method. 
Since the former 
puts restrictions on the ambient variety,
we will only consider examples where the ambient is a projective space.
We are unable to determine precisely
which is the fastest among all the methods
   and which, in the probabilistic case, has highest probability of giving the correct answer.
We just summarize in Table~\ref{table: compare} the running times for some special examples. 
Below is the code from which we obtained the first row of the Table.
{\footnotesize
\begin{Verbatim}[commandchars=&\{\}] 
i27 : loadPackage "&colore{YellowOrange}{CharacteristicClasses}"; loadPackage "&colore{OliveGreen}{FMPIntersectionTheory}"; 
i29 : X = last(P5=ZZ/16411[vars(0..5)],ideal(random(3,P5),random(3,P5),random(4,P5)));
i30 : (time &colore{YellowOrange}{Segre} X,time &colore{blue}{SegreClass} X,time &colore{OliveGreen}{segreClass} X,time &colore{blue}{SegreClass}(X,&colore{blue}{MathMode}=>true));
     &colore{Sepia}{-- used 0.1511 seconds}
     &colore{Sepia}{-- used 1.00936 seconds}
     &colore{Sepia}{-- used 34.1471 seconds}
     &colore{Sepia}{-- used 74.572 seconds}
\end{Verbatim}
} \noindent 

 \begin{table}[htbp]
\centering \footnotesize 
\tabcolsep=1.8pt 
\renewcommand{\arraystretch}{1.1}
\begin{tabular}{|l|rr|rr|}
 \hline 
 \normalsize Input & \begin{tabular}{r}\emph{Characteristic-}\\ \emph{Classes.m2} \end{tabular}& \begin{tabular}{c} \emph{Cremona.m2} \\ (probabilistic) \end{tabular} & \begin{tabular}{r}\emph{FMPIntersection-}\\ \emph{Theory.m2}\end{tabular} &  \begin{tabular}{c} \emph{Cremona.m2} \\ (deterministic) \end{tabular} \\
\hline
\normalsize complete int. of type $(3,3,4)$ in $\PP^5$ &\normalsize 0.15 &\normalsize 1.01 &\normalsize 34.15 &\normalsize 74.57 \\ 
\normalsize rational normal surface $S(1,4)\subset\PP^6_{\mathbb{Q}}$ &\normalsize 1.41 &\normalsize  0.74 &\normalsize 5.32 &\normalsize  0.06 \\ 
\normalsize Grassmannian $\mathbb{G}(1,4)\subset\mathbb{P}^9_{\mathbb{Q}}$ &\normalsize 0.16 &\normalsize 0.09  &\normalsize 0.42 &\normalsize 0.02  \\
\normalsize base locus of $\phi$ in Ex. \ref{cubo} &\normalsize 0.23 &\normalsize 0.44 &\normalsize 6.49 &\normalsize  663.79 \\
\normalsize  $X\subset\PP^{11}$ in Ex. \ref{segreDet} { over $\mathbb{F}_{331^2}$} &\normalsize  65.76  &\normalsize 83.66 &\normalsize  --  &\normalsize  0.85   \\ 
\normalsize  $X\subset\PP^{11}$ in Ex. \ref{segreDet} { over $\mathbb{Z}/16411$} &\normalsize  3.62  &\normalsize 11.61 &\normalsize 198.92   &\normalsize  0.74  \\ 
\hline 
\end{tabular}
 \caption{Run-times to compute Segre classes (all times given in seconds)}
\label{table: compare} 
\end{table}
\end{example}


\providecommand{\bysame}{\leavevmode\hbox to3em{\hrulefill}\thinspace}
\providecommand{\MR}{\relax\ifhmode\unskip\space\fi MR }
\providecommand{\MRhref}[2]{%
  \href{http://www.ams.org/mathscinet-getitem?mr=#1}{#2}
}
\providecommand{\href}[2]{#2}

\end{document}